\documentclass{article}
 \usepackage{arxiv}
\usepackage{float}
\usepackage[utf8]{inputenc} % allow utf-8 input
\usepackage[T1]{fontenc}    % use 8-bit T1 fonts
\usepackage{hyperref}       % hyperlinks
\usepackage{url}            % simple URL typesetting
\usepackage{booktabs}       % professional-quality tables
\usepackage{amsfonts}       % blackboard math symbols
\usepackage{nicefrac}       % compact symbols for 1/2, etc.
\usepackage{microtype}      % microtypography
\usepackage{lipsum}
\usepackage{enumitem}
\usepackage{graphicx}
\usepackage{amsmath}
\usepackage[pagewise]{lineno}
\graphicspath{ {./images/} }
\usepackage{amsthm}
\usepackage[skip=2pt]{caption}
\usepackage[figurename=Fig.]{caption}
 \usepackage[usenames,dvipsnames]{pstricks}
 \usepackage{epsfig}
 \usepackage{pst-grad} % For gradients
 \usepackage{pst-plot} % For axes
 \usepackage{algorithm}
\usepackage[utf8]{inputenc}

\newtheorem{thm}{Theorem}[section]
\newtheorem{lem}{Lemma}[section]
\theoremstyle{definition}

\newtheorem{cor}{Corollary}[section]

\title{Alon-Tarsi Number of Some Regular Graphs}

\author{
  Prajnanaswaroopa S \\
  \texttt{sntrm4@rediffmail.com} \\
\date{Department of Mathematics, Amrita University, Coimbatore-641112}
}

\begin{document}
\maketitle
\begin{abstract}
The Alon-Tarsi number of a polynomial is a parameter related to the exponents of its monomials. For graphs, their Alon-Tarsi number is the Alon-Tarsi number of their graph polynomials. As such, it provides an upper bound on their choice and online choice numbers. In this paper, we obtain the Alon-Tarsi number of some complete multipartite graphs, line graphs of some complete graphs of even order, and line graphs of some other regular graphs.
\end{abstract}
\section{Introduction}
The Combinatorial Nullstellensatz theorem has become a trending method in Algebraic Combinatorics. It is used in a variety of unexpected areas of Combinatorics and Algebra. The theorem, which is an extension or a generalization of the fundamental theorem of algebra to several variables, is quite useful, primarily in graph theory, number theory, and discrete geometry. Here, we use the method (theorem) in graph theory, specifically for graph colorings/list colorings. 

Given a graph $G$ with order and size $n$ and $E$ respectively, we define the average degree of the graph as $a=2\frac{E}{n}$. Note that, for regular graphs, the average degree of the graph equals their maximum and minimum degrees, hence $a=\Delta$ for such graphs, where $\Delta$ is the (maximum) degree of the regular graph. If we give a suitable ordering of the vertices of $G$ in the form $x_1, x_2,\ldots, x_n$, we define the graph polynomial as the product $P=\prod_{i<j}(x_i-x_j)$ where $x_i$ is adjacent to $x_j$. Note that the graph polynomial $P$ is homogenous with degree equal to the size of $G$ and the maximum exponent of any variable equal to the maximum degree of $G$. From the original paper regarding Combinatorial Nullstellensatz by Alon(\cite{ALO}) and the crucial paper by Tarsi and Alon(\cite{ALOT}), we can define the Alon-Tarsi Number of any polynomial  $H=\sum_tc_t\mathbf{y_t}$ with $\mathbf{y_t}=y_{1}^{i_1}y_{2}^{i_2}\ldots y_n^{i_n}$ as $min_{t}(max_{i_k}(y_{1}^{i_1}y_{2}^{i_2}\ldots y_{n}^{i_n}))$, that is, it is the minimum of the highest exponent of any of the monomials in the polynomial. For graphs, the Alon-Tarsi number can be defined as the Alon-Tarsi number of its graph polynomial. (Taking this stance from \cite{ZHU}). From the above min-max relation for the Alon-Tarsi number, we say a monomial with the minimum exponent (exponent refers to the maximum taken over all variables in that monomial) is the Alon-Tarsi monomial. From the references, it is clear that the Alon-Tarsi number of any graph is an upper bound (sometimes strict) on the choice number and hence the chromatic number. The Alon Tarsi number can be derived from the structural properties of the graph as well. Per \cite{ALOT}, it is seen that the Alon Tarsi number is equal to $k+1$, where $k$ is the maximum outdegree of an orientation $D$ of a graph $G$ such that the number of spanning eulerian digraphs of $D$ (digraphs with equal out and in degrees) with even number of edges differ from those spanning eulerian digraphs with odd number of edges. From \cite{CAR} and the primary reference thereof (\cite{SCH}), the Alon-Tarsi number also upper bounds the online list chromatic number or online choice number.

A famous conjecture in the field of list colorings/ list edge colorings is the List Coloring Conjecture, which states that the list edge chromatic number equals the chromatic number for any line graph; or, in other words, the list chromatic index equals the chromatic index for any graph. A graph $G$ is considered $1-$factorizable if its edges can be partitioned into distinct $1-$ factors or perfect matchings. In other words, for a regular graph, $1-$ factorization implies that the number of colors to color the edges equals its degree.

The total graph of a graph $G$, denoted by $T(G)$ \cite{BEH}, is a graph formed by subdividing all the edges of $G$ and connecting the vertices in the subdivided graph that correspond to incident edges of $G$ on the same vertex, as well as vertices which are adjacent in $G$. In this form, it can be seen as the $2-$ distance square of the bipartite graph $S(G)$, the subdivided graph of $G$, with one half square being the line graph $L(G)$ of $G$, and the other half square being $G$ itself. The Total coloring conjecture (TCC) \cite{VIZ}, \cite{BEH1} would mean that $\chi(T(G))\le\Delta(G)+2$. A weaker form of this, the weak TCC \cite{BAS} implies that $\chi(T(G))\le\Delta(G)+3$.\\

One trivial observation from the structure of graph polynomial $P$ is the following.
\begin{lem}
The Alon-Tarsi Number (ATN) of any graph $G$ is at least $\frac{a}{2}$ 
\end{lem}
\begin{proof}
The graph polynomial $P$ of $G$ is homogenous with degree $E$ (size). Therefore, the minimum exponent of any variable will be greater than or equal to $\frac{E}{n}$, where $n$ is the order. By the definition of average degree as given above, the lemma at once follows.
\end{proof}
\section{Theorems}
% \begin{prop}
% Let the vertices of the star graph $G$ with $k+1$ vertices and $k$ edges be labeled $x_1,x_2,\ldots,x_{k+1}$, where $x_1$ is the root vertex, and the rest are pendant vertices. Then, the graph polynomial $G$ has $2^k$ monomials with precisely $\binom{k}{y}$ monomials having an exponent of $x_1$ as $y$.
% \end{prop}
% \begin{proof}
% The proof proceeds by direct multiplication and observation. The graph polynomial is given by $(x_1-x_2)(x_1-x_3)\ldots(x_1-x_{k+1})$. It is clear by observing that any monomial can be formed by first choosing a specified number of $x_1$s and the remaining variables in the factors. Also, each of the monomials so formed is distinct. Hence, their Alon-Tarsi number is $1=1=2$, as the Alon-Tarsi monomial can be chosen as the monomial $x_2x_3\ldots x_{k+1}$.
% \end{proof}
Though the following three results are already implied by the main result in \cite{ALOT}, the approach we use here is relatively straightforward.
\begin{thm}
The Alon-Tarsi monomial of $G=K_{n,n}$ is of the form $c(x_1x_2\ldots x_{2n})^{\frac{n}{2}}\,$ $,c\neq0$ for even $n$. Hence ATN of $G$ is $\frac{n}{2}+1$.
\end{thm}
\begin{proof}
In the graph $G$ we have $a=n\implies \frac{a}{2}=\frac{n}{2}$. Therefore, ATN is bounded below by $\frac{n}{2}+1$. We label all the vertices of one partite set as $x_1,x_2,\ldots,x_n$ and the remaining partite set as $x_{n+1},x_{n+2},\ldots,x_{2n}$. Now, the monomial $c(x_1x_2\ldots x_{2n})^{\frac{n}{2}}$ can be formed by taking all the variables $\frac{n}{2}$ times in the product. Observe that all edges are of the form $(x_i-xj)$, where $i\in\{1,2,\ldots,n\}$ and $j\in\{n+1,n+2,\ldots,2n\}$. The sign of the product of $x_1^{\frac{n}{2}}x_2^{\frac{n}{2}}\ldots x_n^{\frac{n}{2}}$ are all positive as all the individual signs are positive. Now, the signs of $x_{n+1}^{\frac{n}{2}},x_{n+2}^{\frac{n}{2}}\ldots x_{2n}^{\frac{n}{2}}$ can be negative. Nevertheless, as the number of variables in each part (and hence the second part) is even, therefore the sign of the product overall will be $(1)^{n\frac{n}{2}}$, which is positive ($1$). Therefore, the sign of the full product $x_1^{\frac{n}{2}}x_2^{\frac{n}{2}}\ldots x_n^{\frac{n}{2}}x_{n+1}^{\frac{n}{2}},x_{n+2}^{\frac{n}{2}}\ldots x_{2n}^{\frac{n}{2}}$ is positive ($1$). Since this will always be the case, the sum of such monomials will be non-zero. Thus, the Alon-Tarsi monomial of $G$ is of the form $c(x_1x_2\ldots x_{2n})^{\frac{n}{2}}$ for some non-zero $c$. Thus, ATN of $G$ attains its lower bound $\frac{n}{2}+1$. 
\end{proof}
\begin{cor}
The ATN of $G=K_{n,n}$ is $1+\left\lceil\frac{n}{2}\right\rceil.$
\end{cor}
\begin{proof}
The proof is immediate once we note the following three observations:\\
1) The ATN of $K_{n,n}$ for even $n$ is $\frac{n}{2}$.
2) The ATN of $G$ is always $\ge\frac{n}{2}$.
3) $G$ is a subgraph of $K_{n+1,n+1}$.
\end{proof}
\begin{thm}
The ATN of the bipartite graphs $G=K_{m,n}\quad,m<n$ with $n$ even and $(m+n)|mn$ is equal to $\frac{mn}{m+n}+1$. 
\end{thm}
\begin{proof}
In the graph $G$ we have $\frac{a}{2}=\frac{mn}{m+n}$. Therefore, ATN of $G$ is bounded below by $\frac{mn}{m+n}+1$. We label all the vertices of one part as $x_1,x_2,\ldots,x_m$ and the remaining part as $x_{m+1},x_{m+2},\ldots,x_{n}$. The monomial $c(x_1x_2\ldots x_{m+n})^{\frac{mn}{m+n}}$ can be formed by multiplying all the variables $\frac{mn}{m+n}$ times. Observe that the partial product $(x_1x_2\ldots x_{m})^{\frac{mn}{m+n}}$ has a positive sign owing to the positivity of the variables. Whereas, the overall sign of $(x_{m+1}x_{m+2}\ldots x_n)^{\frac{mn}{m+n}}$ is equal to $(-1)^{n\frac{mn}{m+n}}=1$. Therefore, all such monomials are always positive in sign; hence, the sum of these monomials will also be non-zero. Hence, the theorem is implied. 
\end{proof}
\begin{thm}
The ATN of regular bipartite graph $G$ with $2n$ vertices with even $n$ and even degree $\Delta$ is $\frac{\Delta}{2}$.
\end{thm}
\begin{proof}
Again, we label the vertices of $G$ as $x_1,x_2,\ldots, x_n$ for the vertices of one part, and $x_{n+1},x_{n+2},\ldots, x_{2n}$ for the remaining vertices. All edges are of the form $(x_i-x_j)$ where $i$ and $j$ are from different parts. This implies we can form a monomial by choosing each of the variables $\frac{\Delta}{2}$ times. As in the previous theorem, the sign of the product $(x_1x_2\ldots x_n)^{\frac{\Delta}{2}}$ is always positive as all the individual variables are. As for the sign of $(x_{n+1}x_{n+2}\ldots x_{2n})^{\frac{\Delta}{2}}$, we note that as $n$ is even, the combined sign of this partial product would be $(1)^{n\frac{\Delta}{2}}=1$. Therefore, the full sign of the product monomial $(x_1x_2\ldots x_nx_{n+1}x_{n+2}\ldots x_{2n})^{\frac{\Delta}{2}}$ is always positive ($1$). Hence, the sum of such monomials will be non-zero; in equivalent words, ATN of $G$ is $\frac{\Delta}{2}$.
\end{proof}
\begin{thm}
The ATN of the complete $k-$ partite graph $K_{n,n,\ldots(k-times),\ldots,n}$ for even $n$  is equal to $(k-1)\frac{n}{2}$. 
\end{thm}
\begin{proof}
We first orient the graph in an eulerian way (Each vertex has equal out and in degrees), such that for each bipartite graphs induced by any two partite sets of vertices, which are complete, are also oriented in an eulerian way. In this case, the maximum outdegree of the orientation equals half the degree, or $\frac{(k-1)n}{2}$. 

To show that every spanning eulerian subdigraph having even number of edges differs from that having odd number of edges, we proceed as follows.

We see that, in this particular orientation, for spanning subdigraphs, we should have all the vertices of the graph. In addition, the Eulerian-ness forces the subdigraphs to be either to be those subdigraphs formed by concatenation of subdigraphs with even number of edges, or one vertex and an Eulerian subdigraph formed by the other vertices. Thus, the eulerian subdigraphs with even number of edges only surface in the second case, whose number is equal to the order of the graph.

Now, the spanning subdigraphs with even number of edges are:
\begin{enumerate}
\item The empty graph.
\item The graph $G$.
\item Disjoint union of certain eulerian oriented balanced complete multipartite graphs having even number of vertices in each part. 
\item Disjoint union of oriented induced subdigraphs formed by less than $k$ partite sets.
\item Disjoint union of even cycles having length less than or equal to $kn$.
\end{enumerate}
It can be seen, by an induction argument on $k$ with the base case being $k=2$, that number of such subdigraphs easily exceeds $kn$. Thereby, the conditions of Alon-Tarsi theorem are satisfied which instantly gives us the result.
\end{proof}
A related result to the following result was proved using elaborate techniques in \cite{SCHA} for complete graphs of prime degree.
\begin{thm}
The ATN of the line graph of $G=K_n$ for $n=4k\quad,k\in\mathbb{N}$ is $n-1$. Hence, the edge choosability of $K_n$ is $n-1$, or are chromatic edge choosable.
\end{thm}
\begin{proof}
We note that we can factorize the edges of $K_n$ into $\Delta=n-1$ $1-$factors or perfect matchings. We also observe that an edge in any perfect matching is adjacent to exactly two edges in any other perfect matching (the adjacencies are determined by the end vertices). Thus, the line graphs of $ G $ are $ n-1 $ -partite with any two partite sets inducing a regular bipartite graph of degree $ 2 $.

 Let us label the first partite set as $A$ for loss of generality. We can orient $G$ such that any two partite sets induce an eulerian orientation. In addition, this orientation is obtained by oriented the cliques formed by any vertex acyclically, thereby, we have that the edges of the cliques incident at a vertex in $A$ are either oriented towards it or away from it. This can also said to be kernel perfect orientation. In this case, the maximum outdegree of the orientation equals half the degree of $L(G)$, or $\frac{(4k-1)2}{2}=4k-1$. 

To show that every spanning eulerian subdigraph having even number of edges differs from that having odd number of edges, we proceed as follows.

We see that, in this particular orientation, for spanning subdigraphs, we should have all the vertices of the graph. In addition, the Eulerian-ness forces the subdigraphs to be either to be those subdigraphs formed by concatenation of subdigraphs with even number of edges, or one vertex and an Eulerian subdigraph formed by the other vertices. Thus, the eulerian subdigraphs with even number of edges only surface in the second case, whose number is equal to the order of the graph.

Now, the spanning subdigraphs with even number of edges are:
\begin{enumerate}
\item The empty graph.
\item The graph $G$.
\item Disjoint union of certain eulerian oriented sub-line graphs having even number of vertices in each part. 
\item Disjoint union of oriented induced subdigraphs formed by less than $4k$ partite sets.
\item Disjoint union of even cycles having length less than or equal to $4k$.
\end{enumerate}
It can be seen, by an induction argument on $k$  with the base case being $k=1$, that number of such subdigraphs easily exceeds $4k$. Thereby, the conditions of Alon-Tarsi theorem are satisfied which instantly gives us the result.
% We note that we can factorize the edges of $K_n$ into $\Delta=n-1$ $1-$factors or perfect matchings. We label the edges (vertices of line graph) as $x_1,x_2,\ldots,x_{2k}$ for the first perfect matching; $x_{2k+1},x_{2k+2},\ldots,x_{n}$ and so on. We also observe that an edge in any perfect matching is adjacent to exactly two edges in any other perfect matching (the adjacencies are determined by the end vertices). Thus, the line graphs of $ G $ are $ n-1 $ -partite with any two partite sets inducing a regular bipartite graph of degree $ 2 $. Now, from Theorem 2.3, the Alon-Tarsi monomial of the induced graph formed by any two partite sets in the line graph is of the form $c_i(x_{i_1}x_{i_2}\ldots x_{i_{4k}})^{\frac{2}{2}}=c_i(x_{i_1}x_{i_2}\ldots x_{i_{4k}})$ for some non-zero $c_i$ and indices $i_k$s. Since there are $(n-2)$ such regular bipartite graphs (the induced graphs between any two perfect matchings) for each of the partite sets, we get that the final Alon-Tarsi monomial would be of the form $C(x_{i_1}x_{i_2}\ldots x_{i_{2k}})^{n-2}(x_{i_{2k+1}}x_{i_{2k+2}}\ldots x_{4k})^{n-3}\ldots(x_{i_{n-2k+1}}x_{i_{n-2k+2}}\ldots x_{i_{n}})$. Therefore, the Alon-Tarsi number should be $(n-2)+1=n-1$. Hence, the claim follows.
\end{proof}
\begin{thm}
If $G$ is an $n$ order $1-$factorizable regular graph with maximum degree $\Delta$ and $n=4k$ for some integer $k$, then ATN of the line graph of $G$ is equal to $n-1$. Hence, the LCC holds for these graphs.
\end{thm}
\begin{proof}
The proof is similar to the previous theorem. Here, we partition the line graph into $\Delta$ independent sets, each having $\frac{n}{2}=2k$ vertices. Now, the arguments of the previous theorem hold. 

Let us label the first partite set as $A$ for loss of generality. We can orient $G$ such that any two partite sets induce an eulerian orientation. In addition, this orientation is obtained by oriented the cliques formed by any vertex acyclically, thereby, we have that the edges of the cliques incident at a vertex in $A$ are either oriented towards it or away from it. This can also said to be kernel perfect orientation. In this case, the maximum outdegree of the orientation equals half the degree of $L(G)$, or $\frac{(\Delta-1)2}{2}=\Delta-1$. 

To show that every spanning eulerian subdigraph having even number of edges differs from that having odd number of edges, we proceed as follows.

We see that, in this particular orientation, for spanning subdigraphs, we should have all the vertices of the graph. In addition, the Eulerian-ness forces the subdigraphs to be either to be those subdigraphs formed by concatenation of subdigraphs with even number of edges, or one vertex and an Eulerian subdigraph formed by the other vertices. Thus, the eulerian subdigraphs with even number of edges only surface in the second case, whose number is equal to the order of the graph.

Now, the spanning subdigraphs with even number of edges are:
\begin{enumerate}
\item The empty graph.
\item The graph $G$.
\item Disjoint union of certain eulerian oriented sub-line graphs having even number of vertices in each part. 
\item Disjoint union of oriented induced subdigraphs formed by less than $4k$ partite sets.
\item Disjoint union of even cycles having length less than or equal to $4k$.
\end{enumerate}
It can be seen, by an induction argument on $k$  with the base case being $k=1$, that number of such subdigraphs easily exceeds $4k$. Thereby, the conditions of Alon-Tarsi theorem are satisfied which instantly gives us the result.

% and we label the vertices of the line graph as $x_1,x_2,\ldots,x_{2k}$ for the first perfect matching; $x_{2k+1},x_{2k+2},\ldots,x_{n}$ and so on. The edges are now of the form $x_i-x_j$, wghere $i<j$. We also observe that an edge in any perfect matching is adjacent to exactly two edges in any other perfect matching (the adjacencies are determined by the end vertices). Thus, the line graphs of $ G $ are $\Delta$ -partite with any two partite sets inducing a regular bipartite graph of degree $ 2 $. Now, from Theorem 2.3, the Alon-Tarsi monomial of the induced graph formed by any two partite sets in the line graph is of the form $c_i(x_{i_1}x_{i_2}\ldots x_{i_{4k}})^{\frac{2}{2}}=c_i(x_{i_1}x_{i_2}\ldots x_{i_{4k}})$ for some non-zero $c_i$ and indices $i_k$s. Since there are $(\Delta-1)$ such regular bipartite graphs (the induced graphs between any two perfect matchings) for each independent set of the line graph, we get that the final Alon-Tarsi monomial should be of the form $C(x_{i_1}x_{i_2}\ldots x_{i_{2k}})^{\Delta-1}(x_{i_{2k+1}}x_{i_{2k+2}}\ldots x_{4k})^{\Delta-2}\ldots(x_{i_{n-2k+1}}x_{i_{n-2k+2}}\ldots x_{i_{n}})$. This implies that the Alon-Tarsi number should be $(\Delta-1)+1=\Delta$. Hence, $G$ satisfies LCC.
\end{proof}
\begin{cor}
Let $G$ be a graph having order $4k$ that is $1$-factorizable. If $T(G)$ is the total graph of $G$, which is formed by the disjoint union of $G$, line graph $L(G)$, and the edges of the subdivision graph $S(G)$. Then, $ATN(T(G))\le\Delta(G)+2.$ 
\end{cor}
\begin{proof}
The proof is immediate, as the graph $T(G)$ can be seen as the $2$-distance (or square) graph of the subdivision graph, with one half square being $G$ and the other half square being $L(G)$. The interconnecting edges between the two half squares are the edges of the bipartite graph $S(G)$. As such, any minimal Alon-Tarsi monomial will have at most $2$ as an increment in the maximum exponent. Therefore, the Corollary follows.
\end{proof}
\section*{Conflict of Interest}
The author declares no conflict of interest whatsoever.

\end{document}